
  %
  \newcount\fontset
  \fontset=2
  \def\dualfont#1#2#3{\font#1=\ifnum\fontset=1 #2\else#3\fi}
  \dualfont\eightrm {cmr8} {cmr7}
  \dualfont\eightsl {cmsl8} {cmr7}
  \dualfont\eightit {cmti8} {cmti10}
  \dualfont\eightmi {cmmi8} {cmmi10}
  \dualfont\tensc {cmcsc10} {cmcsc10}
  \dualfont\titlefont {cmbx12} {cmbx10}
  \dualfont\eufb {eufb14 scaled 833} {cmsy10}

  \magnification=\magstep1
  \nopagenumbers
  \voffset=2\baselineskip
  \advance\vsize by -\voffset
  \headline{\ifnum\pageno=1 \hfil \else \tensc\hfil
    fell bundles and measurable twisted actions
  \hfil\folio \fi}

  %
  \def\vg#1{\ifx#1\null\null\else \ifx#1,,\else \ifx#1..\else
\ifx#1;;\else \ifx#1::\else \ifx#1''\else \ifx#1--\else
\ifx#1))\else\ifx{#1}{ }{ }\else { }#1\fi\fi\fi\fi\fi\fi\fi\fi\fi}

  \newcount\secno \secno=0
  \newcount\stno 
  \outer\def\section#1{
    \stno=0
    \global\advance\secno by 1
    \vskip0pt plus.1\vsize\penalty-250
    \vskip0pt plus-.1\vsize\bigskip\vskip\parskip
    \message{\number\secno.#1\enspace}
    \noindent{\bf\number\secno.\enspace #1}.}
  \def\state#1 #2\par{\advance\stno by 1\medbreak\noindent
    {\bf\number\secno.\number\stno.\enspace #1.\enspace}{\sl
    #2}\medbreak}
  \def\nstate#1 #2#3\par{\state{#1} {#3}\par
    \edef#2{\number\secno.\number\stno}}
  \def\proof{\medbreak\noindent{\it Proof.\enspace}}
  \def\proofend{\ifmmode\eqno\square\else\hfill\square\medbreak\fi}
  \newcount\itemno \itemno=0
  \def\zitem{\global\advance\itemno by
1\smallskip\item{\ifcase\itemno\or i\or ii\or iii\or iv\or v\or vi\or
vii\or viii\or ix\or x\or xi\or xii\fi)}}
  \def\$#1{#1 $$$$ #1}
  \def\se#1 = #2 ENTAO #3 SENAO #4
FIM{\def\a{#1}\def\b{#2}\ifx\a\b#3\else#4\fi}
  \def\setem#1{\se #1 = {} ENTAO {} SENAO {, #1} FIM}

  %
  \newcount\bibno \bibno=0
  \def\newbib#1{\advance\bibno by 1 \edef#1{\number\bibno}}
  \def\cite#1{{\rm[\bf #1\rm]}}
  \def\scite#1#2{{\rm[\bf #1\rm, #2]}}
  \def\lcite#1{(#1)}
  \def\commonstuff#1{\smallskip \item{[#1]}}
  \def\ATarticle#1{\zarticle #1 xyzzy }
  \def\zarticle#1, author = #2, title = #3, journal = #4, year = #5,
volume = #6, pages = #7, NULL#8 xyzzy {\commonstuff{#1} #2, ``#3'',
\se #7 = {} ENTAO {to appear in {\sl #4}} SENAO {{\sl #4\/} {\bf #6}
(#5), #7.}  FIM}
  \def\ATtechreport#1{\ztechreport #1 xyzzy }
  \def\ztechreport#1, author = #2, title = #3, institution = #4, year
= #5, note = #6, NULL#7 xyzzy {\commonstuff{#1} #2, ``#3''\setem{#6},
#4, #5.}
  \def\ATbook#1{\zbook #1 xyzzy }
  \def\zbook#1, author = #2, title = #3, publisher = #4, year = #5,
volume = #6, series = #7, NULL#8 xyzzy {\commonstuff{#1} #2,
``#3''\setem{#7}\se #6 = {} ENTAO {} SENAO { vol. #6} FIM, #4, #5.}
  \def\ATbooklet#1{\zbooklet #1 xyzzy }
  \def\zbooklet#1, title = #2, author = #3, howpublished = #4, year =
#5, NULL#6 xyzzy {\commonstuff{#1} #3, ``#2'', #4, #5.}

  \def\stress#1{{\it #1}}
  \def\cp{\hbox to 1.8ex{$\times \kern-.45ex\vrule height1.1ex
depth0pt width0.45truept$\hfill}}
  \def\:{\colon}
  \def\*{\otimes}
  \def\x{\times}
  \def\({\bigl(}
  \def\){\bigl)}
  \def\arw{\rightarrow}
  \def\cstar{$C^*$}
  \def\square{\hbox{$\sqcap\!\!\!\!\sqcup$}}
  \def\for#1{,\quad #1}

  \def\goodbreak{\vskip0pt plus.1\vsize\penalty-250 \vskip0pt
plus-.1\vsize\bigskip}
  \def\bundle{Fell bundle\vg}
  \def\bundles{Fell bundles\vg}
  \def\d{\delta}
  \def\D{\partial}
  \def\Bu{\hbox{\eufb B}}
  \def\K{{\cal K}}
  \def\M{{\cal M}}
  \def\UM{{\cal U\kern-0.5pt\M}}

  %
  \newbib{\Banach}
  \newbib{\BS}
  \newbib{\TPA}
  \newbib{\Fell}
  \newbib{\FD}
  \newbib{\PR}

  \null
  \vskip-2\bigskipamount
  \begingroup
  \def\c{\centerline}

  %
  \c{\titlefont CONTINUOUS FELL BUNDLES ASSOCIATED TO}
  \medskip
  \c{\titlefont MEASURABLE TWISTED ACTIONS}\footnote{\null}
    {\eightrm 1991 \eightsl MR Subject Classification:
    \eightrm
    46L05. 
    }

  \baselineskip=10pt
  \eightit

  \bigskip
  \c{\tensc Ruy Exel\footnote{*}{\eightrm Partially supported by CNPq,
Brazil.}}
  \c{Departamento de Matem\'atica}
  \c{Universidade de S\~ao Paulo}
  \c{C.~P.~ 20570}
  \c{01452-990 S\~ao Paulo -- BRAZIL}
  \c{exel@ime.usp.br}

  \bigskip
  \c{\tensc Marcelo Laca\footnote {$\null^\dagger$}{\eightrm Supported
by the Australian Research Council.}}
  \c{Mathematics Department}
  \c{University of Newcastle}
  \c{Newcastle, NSW 2308 -- AUSTRALIA}
  \c{marcelo@math.newcastle.edu.au}

  %
  \bigskip\bigskip\eightrm\baselineskip=3ex
  \midinsert \narrower\narrower
  Given a \underbar{measurable} twisted action of a second-countable,
locally compact group
  {\eightmi G}
  on a separable
  {\eightmi C}*-algebra
  {\eightmi A},
  we prove the existence of a topology on
  {\eightmi A$\x$G}
  making it a \underbar{continuous} \bundle, whose cross sectional
  {\eightmi C}*-algebra is isomorphic to the
Busby--Smith--Packer--Raeburn crossed product.
  \endinsert \endgroup

  \section {Introduction} Let $A$ be a \cstar-algebra and $G$ be a
locally compact group. According to \cite{\BS}, a \stress{twisted
action} of $G$ on $A$ is a pair $(\theta,w)$ of maps
    $
  \theta \: G \arw Aut(A),
    $
  and
    $
  w\: G\x G \arw \UM(A),
    $
  where $Aut(A)$ denotes the automorphism group of $A$ and $\UM(A)$ is
the set of unitary elements in the multiplier algebra $\M(A)$,
satisfying
  \itemno=0
  \zitem $\theta_e$ is the identity automorphism of $A$,
  \zitem $\theta_r(\theta_s(a)) = w(r,s)\theta_{rs}(a)w(r,s)^*$,
  \zitem $w(e,t)=w(t,e)=1$,
  \zitem $\theta_r(w(s,t))w(r,st)=w(r,s)w(rs,t)$,
  \smallskip \noindent for $r,s,t$ in $G$ and $a$ in $A$.

  \goodbreak
  \medskip
  Also following \cite{\BS}, we will say that $(\theta,w)$ is
\stress{measurable} if $\theta$ is {\it strongly measurable} in the
sense that for each $a$ in $A$, the map
    $$
  t \in G \mapsto \theta_t(a) \in A \eqno{(1)}
    $$
  is Borel measurable, and if $w$ is {\it strictly measurable} in the
sense that, for each $a$ in $A$, the maps
    $$
  (r,s) \in G \x G \mapsto aw(r,s) \in A \eqno{(2)}
    $$
    $$
  (r,s) \in G \x G \mapsto w(r,s)a \in A \eqno{(3)}
    $$
  are Borel measurable.

  For convenience, we will restrict ourselves to the treatment of
separable \cstar-algebras and second countable groups, as in
\cite{\BS}, so that, in particular, the various notions of
measurability for $A$ valued maps coincide.

  Finally, $(\theta,w)$ will be termed \stress{continuous} if the maps
(1) -- (3) above are continuous. In this case, it is easy to see that
the conditions of \scite{\TPA}{Definition 3.8} are satisfied, and
hence that we can construct the associated \stress{semi-direct product
bundle} of $A$ by $G$ \scite{\TPA}{Theorem 3.10}.

 \bundles, also frequently referred to as \cstar-algebraic bundles,
were introduced by {J.~M.~G.~Fell}
  \cite{\Fell} (see also \cite{\FD})
  in the 60's.
  Among the many interesting features surrounding this concept, we
would like to point out its relevance for the study of crossed product
\cstar-algebras. In fact \bundles can (and should) be viewed as
intermediate steps in the construction of crossed products. The
procedure being to start by constructing the associated semi-direct
product bundle \scite{\FD}{VIII.4}, \cite{\TPA} and then to consider
its cross sectional algebra \scite{\FD}{VIII.17.2}.

  The available theory of \bundles, of which \cite{\FD} is one of the
most authoritative accounts, does not include, as far as we know, a
systematic study of measurable (as opposed to continuous) bundles.
However, crossed products by measurable twisted actions have been
profitably studied by Packer and Raeburn \cite{\PR}, where they play
an important role in the theory of group actions on
\cstar-algebras. Therefore, it seems plausible that this latter
crossed product construction could be obtained in a similar two step
procedure, involving, as the intermediate step, the construction of a
``measurable'' \bundle.

  Our main point, however, is that, given a measurable action of a
second-countable group on a separable algebra, the associated $L^1$
algebra studied by Busby-Smith \cite{\BS}, as well as the crossed
product of Packer--Raeburn \cite{\PR} can both be obtained from a {\it
continuous} \bundle. This result bears a certain degree of similarity
with the result of S.~Banach, according to which a measurable
homomorphism between complete metric groups is necessarily continuous
\scite{\Banach}{Theorem I.4}.

  Given that the theory of induced representations, as well as the
Mackey normal subgroup analysis has, to a large extent, been
generalized to \bundles\ (see \scite{\FD}{chapters XI and XII}), our
result makes this machinery available to the study of measurable
twisted actions.  This aspect has already been conjectured by Packer
and Raeburn in the introduction to \cite{\PR}.

  The present work was largely developed during the first author's
visit to the Mathematics Department at the University of Newcastle,
whose members he would like to thank for their warm hospitality.

  \section {The Main Results} Let us fix, throughout, a measurable
twisted action $(\theta,w)$ of a locally compact, second-countable
group $G$ on a separable \cstar-algebra $A$. If we let $G^d$ denote
the discrete group obtained by replacing the given topology of $G$ by
the discrete topology, then, following \scite{\TPA}{Theorem 2.7}, we
get a \bundle structure on $A\x G^d$, by means of the following
operations, where we use the notation $a\d_t$ for $(a,t)$ in $A\x
G^d$:
    $$
  (a\d_r)(b\d_s) = a\theta_r(b)w(r,s)\d_{rs} \for a,b\in A,\quad
r,s\in G
    $$
  and
    $$
  (a\d_t)^* = \theta_t^{-1}(a^*)w(t^{-1},t)^*\d_{t^{-1}} \for a\in
A,\quad t\in G^d.
    $$

  We will refer to this bundle as $\Bu(A,G^d,\theta,w)$. Of course it
bears no relationship, whatsoever, with the topological nature of
$G$. Our main result is precisely intended to fill this gap. We will
therefore show that there exists a topology on $A\x G$ for which
$\Bu(A,G^d,\theta,w)$ is a continuous \bundle over $G$. (Whenever $G$
occurs without the superscript `$d$', it is to be thought of as
carrying its originally given topology). In addition, we will show
that the Banach algebra of $L^1$ sections of $\Bu(A,G^d,\theta,w)$ is
isomorphic to $L^1(G,A,\theta,w)$ (see \cite{\BS}) and its cross
sectional \cstar-algebra is isomorphic to the Packer--Raeburn crossed
product $A\cp_{\theta,w}G$, both isomorphisms being canonical.

  One of the main ingredients in the proof of this result is the
``Packer--Raeburn stabilization trick'' \scite{\PR}{Theorem 3.4},
which we briefly describe below, mainly to fix our notation.

  \nstate{Definition} \ExtEq (\scite{\PR}{Definition 3.1}).  The
measurable twisted actions $(\alpha,u)$ and $(\beta,w)$ of $G$ on $A$
are \stress{exterior equivalent} if there exists a strictly measurable
map $v\:G\arw \UM(A)$ such that
  \itemno=0
  \zitem $\beta_t(a) = v_t \alpha_t(a) v_t^* \for a\in A,\quad t\in
G.$
  \zitem $w(r,s) = v_r\alpha_r(v_s)u(r,s)v_{rs}^* \for r,s\in G.$

  \nstate{Proposition} \IsoBundle If $(\alpha,u)$ and $(\beta,w) $are
exterior equivalent then the associated bundles
  $\Bu(A,G^d,\alpha,u)$
  and
  $\Bu(A,G^d,\beta,w)$
  are isomorphic.

  \proof One checks that the map $\phi\:A\x G \arw A\x G$ defined by
$\phi(a,t) = (av_t^*,t)$ is an isomorphism for the respective bundle
structures. In proving this, the identity $\beta_t^{-1}(v_t) =
\alpha_t^{-1}(v_t)$, which follows from \lcite{\ExtEq.i} with $a=
\alpha_t^{-1}(v_t)$, comes in handy.  \proofend

The Packer--Raeburn stabilization trick asserts that, if $(\theta,w)$
is a a measurable twisted action of the second-countable group $G$ on
the separable \cstar-algebra $A$, then there exists a strongly
continuous (untwisted) action $\beta$ of $G$ on $A\*\K$ (where $\K$
denotes the algebra of compact operators on a separable Hilbert
space), such that $(\theta\*1,w\*1)$ is exterior equivalent to
$(\beta,1)$.  Therefore, by \lcite{\IsoBundle}, the bundles
    $$
  \Bu_0 : = \Bu(A\*\K,G^d,\theta\*1,w\*1)
    $$
  and
    $$
  \Bu_1 : = \Bu(A\*\K,G^d,\beta,1)
    $$
  are isomorphic. Now, since $\beta$ is strongly continuous, the
product topology on $A\x G$\/ makes $\Bu_1$ a continuous \bundle over
$G$ \scite{\TPA}{Theorem 3.10}. Thus we can make $\Bu_0$ into a
continuous \bundle over $G$ by transferring the topology from $\Bu_1$
to $\Bu_0$ via the isomorphism given by \lcite{\IsoBundle}.

 On the other hand, if $p$ denotes a minimal projection in $\K$, one
has that $\Bu(A,G^d,\theta,w)$ sits naturally in $\Bu_0$ as the
sub-bundle $(A\*p)\x G^d$.

  The crucial point in our argument is to show that this sub-bundle,
with the inherited topology, is a continuous \bundle.

  In order to verify this claim, let us fix some notation. First of
all, given the two distinct bundle structures on $(A\*\K)\x G^d$, let
us agree to denote the element $(x,t)$ in $(A\*\K)\x G^d$ by $x\D_t$
when it is viewed as an element of $\Bu_1$ while retaing the notation
$x\d_t$ when $\Bu_0$ is concerned. Secondly, let us denote by $v$ a
given Borel map
    $
  v\:G \arw \UM(K\*A)
    $
 implementing the equivalence \lcite{\ExtEq} between
$(\theta\*1,w\*1)$ and $(\beta,1)$, which exists by the stabilization
trick.

  By Proposition \lcite{\IsoBundle}, the map
    $$
  \phi(x\d_t) = x v_t^*\D_t \for x\in A\*\K,\quad t\in G
    $$
  is an isomorphism from $\Bu_0$ to $\Bu_1$.

  Recall that any multiplier of the unit fiber algebra of a \bundle
extends to a multiplier of order $e$ (where $e$ denotes the group
unit) of the bundle concerned \scite{\FD}{VIII.3.8}. In particular,
$1\*p$, viewed as a multiplier of $A\*\K$, extends to a multiplier of
$\Bu_1$, which we will denote by $\pi$.

  \state {Lemma} For each $t$ in $G$ one has
    $$
 \phi\((A\*p)\d_t\) = \pi \( (A \* \K) \D_t \)\pi.
    $$

  \proof We have
    $$
  \pi \( (A \* \K) \D_t \)\pi
  =
  (1\*p) \D_e (A\*\K)\D_t (1\*p)\D_e
  \$=
  (1\*p) (A\*\K)\beta_t(1\*p)\D_t
  =
  (1\*p) (A\*\K)v_t(\theta_t\* 1)(1\*p)v_t^*\D_t
  \$=
  (1\*p) (A\*\K)v_t(1\*p)v_t^*\D_t
  =
  (A\*p) v_t^*\D_t
  =
  \phi\bigl((A\*p)\d_t\bigr).
  \proofend
    $$

  This brings us to our first main result.

  \nstate{Theorem} \Topology Let $(\theta,w)$ be a measurable twisted
action of the locally compact, second-countable group $G$ on the
separable \cstar-algebra $A$. Also let $\Bu_1$ be as above. Then
$\Bu(A,G^d,\theta,w)$, equipped with the topology induced by the
embedding
    $$
  \Bu(A,G^d,\theta,w) \arw \Bu_1
    $$
  given by
    $$
   a\d_t \mapsto (a\*p)v_t^*\D_t,
    $$
  is a continuous \bundle bundle over $G$.

  \proof The proof consists in showing that the collection of
subspaces
    $$
  (A\*p)v_t^*\D_t \subseteq (A\*\K)\D_t \for t\in G,
    $$
  forms a continuous bundle of Banach spaces over $G$
\scite{\FD}{II.13.1}, in the sense that one can find a continuous
section $\gamma$ passing through any preassigned element
$(a\*p)v_{t_0}^*\D_{t_0}$, and such that $\gamma(t) \in
(A\*p)v_t^*\D_t$ for all $t$ in $G$ \scite{\FD}{II.13.18},
\scite{\TPA}{Proposition 3.3}.
   For this, it suffices to take a section $\sigma$ of $\Bu_1$ such
that
  $\sigma({t_0}) = (a\*p)v_{t_0}^*\D_{t_0}$.
  Since the left and right actions of multipliers are continuous maps
on the bundle \scite{\FD}{VIII.2.14}),
  $\gamma(t) = \pi\sigma(t)\pi$ gives the desired section. \proofend

  Let us denote by $\Bu(A,G,\theta,w)$ (omitting the superscript in
$G^d$), the continuous bundle over $G$, arising from
\lcite{\Topology}.

  \state{Theorem} Let $A$, $G$, $\theta$ and $w$ be as in \lcite
{\Topology}. Then the formula $\psi(f)(t)=f(t)\delta_t$, for $f$ in
$L^1(G,A)$ and $t$ in $G$, gives a Banach *-algebra isomorphism
    $$
  \psi\: L^1(G,A,\theta,w) \arw L^1\(\Bu(A,G,\theta,w)\),
    $$
  which, in turn, induces an isomorphism between the crossed product
$A\cp_{\theta,w}G$ and the cross sectional \cstar-algebra of
$\Bu(A,G,\theta,w)$.

  \proof Under the identification provided by Theorem
\lcite{\Topology} we may write $\psi(f)(t) = (f(t)\*p)v_t^*\D_t$.
Observe that, if $f$ is in $L^1(G,A,\theta,w)$ then $\psi(f)$ is an
integrable section of $\Bu(A,G,\theta,w)$ because $v$ is strictly
Borel measurable.  The same reasoning applies to prove the converse
and hence $\psi$ gives an isometric linear isomorphism
    $$
  \psi\:L^1(G,A,\theta,w) \arw L^1\bigl( \Bu(A,G,\theta,w) \bigr).
    $$ It is now easy to show that $\psi$ is also a Banach *-algebra
isomorphism.

The description of $A\cp_{\theta,w}G$ which better suits our purposes
is that given in \scite{\PR}{Remark 2.6}, where $A\cp_{\theta,w}G$ is
described as the enveloping \cstar-algebra of
$L^1(G,A,\theta,w)$. Since, on the other hand, the cross sectional
\cstar-algebra of $\Bu(A,G,\theta,w)$ is the enveloping \cstar-algebra
of the algebra of integrable sections of this bundle, we see that the
last part of the statement follows by taking the enveloping
\cstar-algebras of the corresponding Banach algebras.  \proofend

  \bigbreak
  \centerline{\tensc References}
  \nobreak\medskip
  \frenchspacing

\ATbook{\Banach,
  author = {S. Banach},
  title = {Th\'eorie des Op\'erations Lin\'eaires},
  publisher = {Hafner Publishing Co., New York},
  year = {1932},
  volume = {},
  series = {},
  NULL = {},
  }

\ATarticle{\BS,
  author = {R. C. Busby and H. A. Smith},
  title = {Representations of Twisted Group Algebras},
  journal = {Trans. Amer. Math. Soc.},
  year = {1970},
  volume = {149},
  pages = {503--537},
  NULL = {},
  }

\ATtechreport{\TPA,
  author = {R. Exel},
  title = {Twisted Partial Actions, A Classification of Stable
$C^*$-Algebraic Bundles},
  institution = {Universidade de S\~ao Paulo},
  year = {1994},
  note = {preprint},
  NULL = {},
  atrib = {N},
  }

\ATbooklet{\Fell,
  title = {An extension of Mackey's method to Banach *-Algebraic
Bundles},
  author = {J. M. G. Fell},
  howpublished = {Memoirs Amer. Math. Soc. vol. 90},
  year = {1969},
  NULL = {},
  }

\ATbook{\FD,
  author = {J. M. G. Fell and R. S. Doran},
  title = {Representations of *-algebras, locally compact groups, and
Banach *-algebraic bundles},
  publisher = {Academic Press},
  year = {1988},
  volume = {125 and 126},
  series = {Pure and Applied Mathematics series},
  NULL = {},
  }

\ATarticle{\PR,
  author = {J. A. Packer and I. Raeburn},
  title = {Twisted crossed products of $C^*$-algebras},
  journal = {Math. Proc. Cambridge Philos. Soc.},
  year = {1989},
  volume = {106},
  pages = {293--311},
  NULL = {},
  }

  \vskip 2cm
  \rightline{June 1995}

  \bye